\documentclass{amsart}
\usepackage{amsfonts}

\usepackage{graphicx}
\usepackage{amscd}

\newtheorem{theorem}{Theorem}
\theoremstyle{plain}
\newtheorem{acknowledgement}{Acknowledgement}

\newtheorem{definition}{Definition}

\newtheorem{remark}{Remark}

\numberwithin{equation}{section}

\input{tcilatex}

\begin{document}
\title[On two variable Drichlet's $q$-$L$-series]{On the two-variable
Drichlet $q$-$L$-series}
\author{YILMAZ\ SIMSEK}
\curraddr{YILMAZ\ SIMSEK\\
Mersin University, Faculty of Science, Department of Mathematics 33343
Mersin, Turkey}
\email{ysimsek@mersin.edu.tr}
\author{DAEYEOUL KIM}
\address{DAEYEOUL KIM\\
KAIST, Department of Mathematics, Tajeon 305-701, S. Korea }
\email{dykim@math.chonbukac.kr}
\author{ SEOG-HOON RIM}
\address{ SEOG-HOON RIM\\
Department of Mathematics Education, Kyungpook National University, Taegu
702-701, S. Korea}
\email{e-mail:shrim@knu.ac.kr}

\begin{abstract}
In this study, we construct the two-variable Dirichlet $q$-$L$-function and
the two-variable multiple Dirichlet-type Changhee $q$-$L$-function. These
functions interpolate the $q$-Bernoulli polynomials and generalized Changhee 
$q$-Bernoulli polynomials. By using the Mellin transformation, we give an
integral representation for the two-variable multiple Dirichlet-type
Changhee $q$-$L$-function. We also obtain relations between the Barnes-type $%
q$-zeta function and the two-variable multiple Dirichlet-type Changhee $q$-$%
L $-function.
\end{abstract}

\maketitle

\section{Introduction}

In his paper \cite{Barnes}, Barnes defined multiple zeta function. Barnes'
multiple zeta function $\zeta _{r}(s,w\mid a_{1},...,a_{r})$ depends on
parameters $a_{1},...,a_{r}$ that will be taken positive throughout this
paper. It is defined by the series 
\begin{equation}
\zeta _{r}(s,w\mid a_{1},...,a_{r})=\sum_{m_{1},...,m_{r}}^{\infty
}(w+m_{1}a_{1}+...+m_{r}a_{r})^{-s},  \label{Eq-12}
\end{equation}%
where $\func{Re}(w)>0,\func{Re}(s)>r$.

The Barnes' multiple Bernoulli polynomials $B_{n}(x,r\mid a_{1},...,a_{r})$,
cf. \cite{Barnes}, are defined by 
\begin{equation}
\frac{t^{r}e^{xt}}{\prod_{j=1}^{r}(e^{a_{j}t}-1)}=\sum_{n=0}^{\infty
}B_{n}(x,r\mid a_{1},...,a_{r})\frac{t^{n}}{n!},  \label{Eq-13}
\end{equation}
for $\mid t\mid <1$.

By (\ref{Eq-12}) and (\ref{Eq-13}), it is easy to see that 
\begin{equation*}
\zeta _{r}(-m,w\mid a_{1},...,a_{r})=\frac{(-1)^{r}m!}{(r+m)!}%
B_{r+m}(w,r\mid a_{1},...,a_{r}),
\end{equation*}%
for $w>0$ and $m$ is a positive integer ( for detail see \cite{Cherednik}, %
\cite{E. Friedman and S. Ruijsenaars}, \cite{S. N. M. Ruijsenaars}, \cite{H.
M. Srivastava and P. W. Karlsson}, \cite{Kim9}, \cite{Kim -simsek }, \cite%
{Kim-simsek2}, \cite{Kim-simsek-sri}).

Recently, many mathematicians and physicians have studied on zeta functions,
multiple zeta functions, multiple $L$-series and multiple Bernoulli numbers
and polynomials due to their importance. These functions and numbers are
used in Number Theory, Complex Analysis and Mathematical Physics, $p$-adic
Analysis and other areas (for detail see \cite{Barnes}, \cite{E. T. Wittaker
and G. N. Watson}, \cite{Apostol}, \cite{Washington}, \cite{Cherednik}, \cite%
{Kim5}, \cite{Kim7}, \cite{Kim9}, \cite{T. M. Rassia and H. M. Srivastava}, %
\cite{Y. Simsek}, \cite{C. A. Nelson and M. G. Gartley-1}, \cite{P. T. Young}
).

In \cite{K. Matsumoto}, Matsumoto studied on general multiple zeta functions
of several variables, involving both Barnes multiple zeta functions and
Euler-Zagier sums as special cases.

In his paper, Ota\cite{K. Ota} gave Kummer-type congruences for derivatives
of Barnes' multiple Bernoulli polynomials. He also generalized these
congruences to derivatives of Barnes' multiple Bernoulli polynomials by an
elementary method and gave a $p$-adic interpolation of them.

By using non-Archimedean $q$-integration, Kim \cite{Kim7} defined multiple
Changhee $q$-Bernoulli polynomials which form a $q$-extension of Barnes'
multiple Bernoulli polynomials. He also constructed the Changhee $q$-zeta
functions (which gives $q$-analogues of Barnes' multiple zeta functions). He
found relations between the Changhee $q$-zeta function and Daehee $q$-zeta
function.

In \cite{P. T. Young}, Young gave some $p$-adic integral representation for
the two-variable $p$-adic $L$-functions. For powers of the Teichm\"{u}ller
character, he used the integral representation to extend the $L$-function to
the large domain, in which it is a meromorphic function in the first
variable and an analytic element in the second. These integral
representations imply systems of congruences for the generalized Bernoulli
polynomials.

In \cite{Kim-simsek-sri}, by using $q$-Volkenborn integration and uniform
differentiable on $\mathbb{Z}_{p}$, Kim, Simsek and Srivastava constructed $%
p $-adic $q$-zeta functions. These functions interpolate the $q$-Bernoulli
numbers and polynomials. The value of $p$-adic $q$-zeta functions at
negative integers were given explicitly. They also defined new generating
functions of $q$-Bernoulli numbers and polynomials. By using these
functions, they proved analytic continuation of $q$-$L$-series. These
generating functions also interpolate the Barnes-type Changhee $q$-Bernoulli
numbers with attached to Dirichlet character as well. By applying the Mellin
transformation, they obtained relations between the Barnes-type $q$-zeta
function and new Barnes-type Changhee $q$-Bernolli numbers.

In \cite{Kim-Rim1}, Kim and Rim constructed two-variable $L$-functoin, $%
L(s,x\mid \chi )$. They showed that this function interpolates the
generalized Bernoulli polynomials associated with $\chi $. By the Mellin
transforms, they gave the complex integral representation for the
two-variable Dirichlet $L$-function. They also found some properties of the
two-variable Dirichlet $L$-function.

In \cite{Kim-simsek2}, Kim constructed the two-variable $p$-adic $q$-$L$%
-function which interpolates the generalized $q$-Bernoulli polynomials
associated with Dirichlet character. He also gave some $p$-adic integrals
representation for this two-variable $p$-adic $q$-$L$-function and derived $%
q $-extension of the generalized formula of Diamond and Ferro and Greenberg
for the two variable $p$-adic $L$-function in terms of the $p$-adic gamma
and log gamma function.

In this paper, more precisely, we define and prove the following results:

We define the two-variable Dirichlet-type Changhee $q$-$L$-function as
follows:

\begin{definition}
Let $\chi $ be a Dirichlet character of conductor $f\in \mathbb{Z}^{+}$. 
\begin{equation*}
L_{q}(s,x\mid \chi ;w_{1})=w_{1}\sum_{n=0}^{\infty }\frac{\chi (n)q^{w_{1}n}%
}{[w_{1}n]^{s}}.
\end{equation*}
\end{definition}

\begin{theorem}
Let $s,w_{1}\in \mathbb{C},$ with $\func{Re}(w_{1})>0.$%
\begin{equation*}
L_{q,r}(s,x\mid \chi ;w_{1})=[f]^{-s}\sum_{a=1}^{f}\chi (a)\zeta _{q^{f}}(s,%
\frac{x+w_{1}a}{f}\mid w_{1}),
\end{equation*}%
where $\zeta _{q}(s,x\mid w_{1})$ is denoted by the Barnes-type Changhee $q$%
-zeta function, which is defined by (see \cite{Kim4}, \cite{Kim5}, \cite%
{Kim7}, \cite{Kim9}, \cite{Kim15}, \cite{Kim-simsek-sri}): For $s\in \mathbb{%
C}$, we have 
\begin{equation*}
\zeta _{q}(s,w\mid w_{1})=-\frac{(1-q)^{s}}{s-1}\frac{1}{\log q}%
+w_{1}\sum_{n=0}^{\infty }\frac{q^{w_{1}n+w}}{[w_{1}n+w]^{s}}.
\end{equation*}
\end{theorem}

A relationship between $L_{q}(s,x\mid \chi ;w_{1})$ and generalized Changhee 
$q$-Bernoulli numbers are given as follows:

\begin{theorem}
Let $\chi $ be a Dirichlet character of conductor $f\in \mathbb{Z}^{+},$ let 
$s,w_{1}\in \mathbb{C},$ with $\func{Re}(w_{1})>0$ and let $n\in \mathbb{Z}%
^{+}$. Then 
\begin{equation*}
L_{q}(1-n,x\mid \chi ;w_{1})=-\frac{\beta _{n,\chi }(x:q\mid w_{1})}{n},
\end{equation*}%
where $\beta _{n,\chi }(x,q\mid w_{1})$ is given by the following generating
function%
\begin{equation*}
-tw_{1}\sum_{n=1}^{\infty }\chi
(n)q^{nw_{1}}e^{[nw_{1}+x]t}=\sum_{n=1}^{\infty }\frac{\beta _{n}(x:q\mid
w_{1})}{n!}t^{n}
\end{equation*}%
(For these polynomials see\ \cite{Kim4}, \cite{Kim5}, \cite{Kim7}, \cite%
{Kim9}, \cite{Kim15}, \cite{Kim-simsek-sri}).
\end{theorem}

The two-variable Dirichlet-type multiple Changhee $q$-$L$-functions are
defined as follows.

\begin{definition}
Let $s,w_{1},...,w_{r}\in \mathbb{C},$ with $\func{Re}(w_{j})>0$, $%
j=1,2,...,r$. For a Dirichlet character $\chi $ with conductor $f\in \mathbb{%
Z}^{+}$, we define 
\begin{equation*}
L_{q,r}(s,x\mid \chi ;w_{1},...,w_{r})=\left( \prod_{j=1}^{r}w_{j}\right)
\sum_{n_{1},n_{2},...,n_{r}=1}^{\infty }\frac{\left( \prod_{k=1}^{r}\chi
(n_{k})\right) q^{\left( \sum_{m=1}^{r}w_{m}n_{m}\right) }}{%
^{[x+\sum_{m=1}^{r}w_{m}n_{m}]^{s}}}.
\end{equation*}
\end{definition}

\begin{theorem}
Let $\chi $ be a Dirichlet character of conductor $f\in \mathbb{Z}^{+}$. Let 
$s,w_{1},...,w_{r}\in \mathbb{C},$ with $\func{Re}(w_{j})>0$, $j=1,2,...,r$.
Then%
\begin{eqnarray*}
L_{q,r}(s,x &\mid &\chi
;w_{1},...,w_{r})=[f]^{r-s}\sum_{a_{1},...,a_{r}=1}^{f}\left(
\prod_{k=1}^{r}\chi (a_{k})\right) \\
\times \zeta _{q^{f},r}(s,\frac{x+w_{1}a_{1}+...+w_{r}a_{r}}{f} &\mid
&w_{1},...,w_{r}),
\end{eqnarray*}%
where $\zeta _{q,r}(s,w\mid w_{1},w_{2},...,w_{r}),$ the Barnes-type
multiple Changhee $q$-zeta functions, is defined by 
\begin{equation*}
\zeta _{q,r}(s,w\mid
w_{1},w_{2},...,w_{r})=\sum_{n_{1},n_{2},...,n_{r}=0}^{\infty }\frac{%
q^{w+n_{1}w_{1}+n_{2}w_{2}+...+n_{r}w_{r}}}{%
^{[w+n_{1}w_{1}+n_{2}w_{2}+...+n_{r}w_{r}]^{s}}},\text{ }
\end{equation*}%
$\Re (w)>0,$ $q\in C$ \ with $\mid q\mid <1$ (For the Barnes-type Changhee
multiple $q$-zeta functions see\ \cite{Kim4}, \cite{Kim5}, \cite{Kim7}, \cite%
{Kim9}, \cite{Kim15}, \cite{Kim-simsek-sri}).
\end{theorem}

We note that $\zeta _{q,r}(s,w\mid w_{1},w_{2},...,w_{r})$\ is analytic
continuation for $\func{Re}(s)>r.$

The numbers $L_{q,r}(-n,\chi \mid w_{1},...,w_{r})$, ( $n>0$ ) are given
explicitly by Theorem below.

\begin{theorem}
Let $\chi $ be a Dirichlet character of conductor $f\in \mathbb{Z}^{+}$. Let 
$w_{1},...,w_{r}\in \mathbb{C},$ with $\func{Re}(w_{j})>0$, $j=1,2,...,r$
and let $n\in \mathbb{Z}^{+}$. Then%
\begin{equation*}
L_{q,r}(-n,x\mid \chi ;w_{1},...,w_{r})=(-1)^{r}\frac{n!}{(n+r)!}B_{n,\chi
}^{(r)}(x:q\mid w_{1},...,w_{r}),
\end{equation*}%
where $B_{n,\chi }^{(r)}(x,q\mid w_{1},...,w_{r}),$ the Barnes-type multiple
Changhee $q$-Bernoulli polynomials, is defined by the following generating
function%
\begin{eqnarray*}
&&(-t)^{r}(\prod_{i=1}^{r}w_{i})\sum_{n_{1},n_{2},...,n_{r}=0}^{\infty
}q^{x+n_{1}w_{1}+n_{2}w_{2}+...+n_{r}w_{r}}e^{[x+n_{1}w_{1}+n_{2}w_{2}+...+n_{r}w_{r}]t}
\\
&=&\sum_{n=0}^{\infty }\frac{B_{n}^{(r)}(x:q\mid w_{1},w_{2},...,w_{r})t^{n}%
}{n!}\text{ \ (}\mid t\mid <2\pi \text{ )}
\end{eqnarray*}%
(For $B_{n,\chi }^{(r)}(x,q\mid w_{1},...,w_{r})$ see\ \cite{Kim4}, \cite%
{Kim5}, \cite{Kim7}, \cite{Kim9}, \cite{Kim15}, \cite{Kim-simsek-sri}).
\end{theorem}

\section{Definition and Notations}

Let $\mathbb{C}$ be the set of complex numbers and $z\in \mathbb{C}$ $.$The
classical Bernoulli polynomials $B_{n}(z)$ are defined by means of of the
generating function (for detail see \cite{E. T. Wittaker and G. N. Watson}, %
\cite{K. Shiratani and S. Yamamoto}, \cite{Apostol}, and \cite{Kim18}, \cite%
{Kim-simsek-sri}): 
\begin{equation}
F(t,x)=\frac{te^{zt}}{e^{t}-1}=\sum_{n=0}^{\infty }B_{n}(z)\frac{t^{n}}{n!}%
\text{ \ (}\mid t\mid <2\pi \text{ ).}  \label{eq-1}
\end{equation}%
Putting $z=0$ into (\ref{eq-1}), $B_{n}(0)=B_{n}$ is the usual $n$th
classical Bernoulli number. The classical Bernoulli numbers are defined by
means of the generating function: 
\begin{equation}
F(t)=\frac{t}{e^{t}-1}=\sum_{n=0}^{\infty }B_{n}\frac{t^{n}}{n!}\text{ \ (}%
\mid t\mid <2\pi \text{ ).}  \label{q-2}
\end{equation}

In \cite{Carlitz}, Carlitz defined $q$-extensions of these classical
Bernoulli numbers and polynomials. Carlitz's $q$-Bernoulli numbers, $\beta
_{n}=\beta _{n}(q)$ are defined by\cite{Carlitz} 
\begin{equation*}
\beta _{0}=1,q(q\beta +1)^{n}-\beta _{n}=\left\{ 
\begin{array}{c}
1,\text{ if }n=1 \\ 
0,\text{ if }n>1,%
\end{array}%
\right.
\end{equation*}%
with the usual convention about replacing $\beta ^{n}$ by $\beta _{n}$.

Carlitz's $q$-Bernoulli polynomials $\beta _{n}(x:q)$ are defined as follows
(\cite{Kim4}, \cite{Kim7}, \cite{Kim-simsek-sri}) 
\begin{equation*}
\beta _{n}(x:q)=\sum_{k=0}^{n}\left( 
\begin{array}{c}
n \\ 
k%
\end{array}%
\right) \beta _{k}q^{kx}[x]^{n-k}.
\end{equation*}%
Thus we note that%
\begin{equation*}
\lim_{q\rightarrow 1}\beta _{n}(q)=B_{n},
\end{equation*}%
and 
\begin{equation*}
\lim_{q\rightarrow 1}\beta _{n}(x:q)=B_{n}(x).
\end{equation*}

Let $\mathbb{Z}^{+}$\ be the set of positive integer numbers. Let $\chi $ be
a Dirichlet character of conductor $f\in \mathbb{Z}^{+}$. Then the
generalized Bernoulli numbers $B_{n,\chi }$ are defined by means of the
generating function%
\begin{equation}
F_{\chi }(t)=\sum_{a=0}^{f-1}\frac{\chi (a)te^{at}}{e^{ft}-1}%
=\sum_{n=0}^{\infty }B_{n,\chi }\frac{t^{n}}{n!}\text{ \ (}\mid t\mid <2\pi 
\text{ ).}  \label{eq-1i}
\end{equation}%
The generalized Bernoulli polynomials $B_{n,\chi }(x)$ are defined by means
of the generating function%
\begin{equation}
F_{\chi }(t,x)=\sum_{a=0}^{f-1}\frac{\chi (a)te^{(a+x)t}}{e^{ft}-1}%
=\sum_{n=0}^{\infty }B_{n,\chi }(x)\frac{t^{n}}{n!}\text{ \ (}\mid t\mid
<2\pi \text{ ).}  \label{q-1ii}
\end{equation}%
Thus these polynomials, $B_{n,\chi }(x),$ are defined as follows: 
\begin{equation}
B_{n,\chi }(x)=\sum_{k=0}^{n}\left( 
\begin{array}{c}
n \\ 
k%
\end{array}%
\right) B_{k,\chi }x^{n-k}\text{.}  \label{eq-3}
\end{equation}

The Riemann zeta function is defined by the series%
\begin{equation*}
\zeta (s)=\sum_{n=0}^{\infty }\frac{1}{n^{s}},
\end{equation*}%
where $s\in \mathbb{C}$ with $\func{Re}(s)>0.$

The Hurwitz zeta function is defined by the series%
\begin{equation*}
\zeta (s,x)=\sum_{n=1}^{\infty }\frac{1}{(n+x)^{s}},
\end{equation*}%
where $s\in \mathbb{C}$ with $\func{Re}(s)>0$ (\cite{Apostol}, \cite{E. T.
Wittaker and G. N. Watson}, \cite{Washington}, \cite{Kim-simsek-sri}, \cite%
{Kim-Rim1}). We note that $\zeta (s,1)=\zeta (s)$. These functions are
analytic continuation for $\func{Re}(s)>1$. Relation between $\zeta (s,x)$
and Bernoulli polynomials, $B_{n}(x)$ is given as follows:

For $n\in \mathbb{Z}^{+},$%
\begin{equation*}
\zeta (1-n,x)=-\frac{B_{n}(x)}{n}
\end{equation*}%
(for detail see \cite{Apostol}, \cite{E. T. Wittaker and G. N. Watson}, \cite%
{Kim18}).

Let $\chi $ be a Dirichlet character of conductor $f\in \mathbb{Z}^{+}$. The
Dirichlet $L$-function is defined by%
\begin{equation*}
L(s,\chi )=\sum_{n=0}^{\infty }\frac{\chi (n)}{n^{s}}.
\end{equation*}%
Relation between $L(s,\chi )$ and the generalized Bernoulli numbers, $%
B_{n,\chi }$ is given as follows: For $n\in \mathbb{Z}^{+},$%
\begin{equation*}
L(1-n,\chi )=-\frac{B_{n,\chi }}{n}.
\end{equation*}

The two-variable $L$-functoin is defined as follows:

\begin{definition}
(\cite{Kim-Rim1}) Let $\chi $ be a Dirichlet character of conductor $f\in 
\mathbb{Z}^{+}$ and $s\in \mathbb{C.}$%
\begin{equation*}
L(s,x\mid \chi )=\sum_{n=0}^{\infty }\frac{\chi (n)}{(n+x)^{s}}.
\end{equation*}
\end{definition}

We note that $L(s,1\mid \chi )=L(s,\chi )$.

$L(s,x\mid \chi )$ is analytic continuation\ in $\mathbb{C}$\ with only
simple pole at $s=1$ (\cite{Kim-Rim1}). Relation between $L(s,x\mid \chi )$
and the generalized Bernoulli numbers, $B_{n,\chi }$ is given as follows:
for $n\in \mathbb{Z}^{+},$%
\begin{equation*}
L(1-n,x\mid \chi )=-\frac{B_{n,\chi }}{n}
\end{equation*}%
( For detail see \cite{Kim-Rim1}).

A sequence of $p$-adic rational numbers as multiple Changhee $q$-Bernoulli
numbers and polynomials are defined as follows\cite{Kim7}, \cite{Kim14}:

Let $a_{1},...,a_{r}$ be nonzero elements of the $p$-adic number field and
let $z\in \mathbb{C}_{p}$. 
\begin{equation}
\beta _{n}^{(r)}(w:q\mid a_{1},...,a_{r})=\frac{1}{\prod_{j=1}^{r}a_{j}}%
\int_{\mathbb{Z}_{p}^{r}}[w+\sum_{j=1}^{r}a_{j}x_{j}]^{n}d\mu _{q}(x),
\label{Eq-14}
\end{equation}%
where 
\begin{equation*}
\int_{\mathbb{Z}_{p}^{r}}f(x)d\mu _{q}(x)=\int_{\mathbb{Z}_{p}}\int_{\mathbb{%
Z}_{p}}...\int_{\mathbb{Z}_{p}}f(x)d\mu _{q}(x_{1})d\mu _{q}(x_{2})...d\mu
_{q}(x_{r})
\end{equation*}%
( see \cite{Kim4}, \cite{Kim5}, \cite{Kim7}, \cite{Kim9}, \cite{Kim14},\cite%
{Kim15}). It is easily observed from (\ref{Eq-14}) that 
\begin{eqnarray}
\beta _{n}^{(r)}(w &:&q\mid a_{1},...,a_{r})=\frac{1}{(1-q)^{n}}%
\sum_{l=0}^{n}\left( 
\begin{array}{c}
n \\ 
l%
\end{array}%
\right) (-1)^{l}q^{wl}\prod_{j=1}^{r}\frac{(l+\frac{1}{a_{j}})}{[la_{j}+1]}
\label{Eq-15} \\
&=&\sum_{l=0}^{n}\left( 
\begin{array}{c}
n \\ 
l%
\end{array}%
\right) [w]^{n-l}q^{wl}\beta _{l}^{(r)}(q\mid a_{1},...,a_{r}),\text{ for }%
n\in \mathbb{Z}^{+}.  \notag
\end{eqnarray}

By using (\ref{Eq-13}) and (\ref{Eq-15}), we note that 
\begin{equation*}
\lim_{n\rightarrow \infty }\beta _{n}^{(r)}(w:q\mid
a_{1},...,a_{r})=B_{n}(w,r\mid a_{1},...,a_{r}).
\end{equation*}%
\begin{equation*}
\beta _{n}^{(r)}(w:q\mid 1,1,...,1)=\beta _{n}^{(r)}(w:q),
\end{equation*}%
where $\beta _{n}^{(r)}(w:q)$ are the $q$-Bernoulli polynomials of order $r$
(see\cite{Kim5}), which are reduced to the ordinary Bernoulli polynomials of
higher order $B_{n}^{(r)}(w)$ if $q=1$ (see for detail \cite{Kim8}, \cite%
{Kim13}, \cite{Kim14}).

In his paper \cite{Kim7}, Kim defined $q$-version of each of the functions $%
F(t)$ and $F(t,x)$ occurring in (\ref{eq-1}) and (\ref{q-2}), respectively.
These generating functions, $F_{q}(t)$ of $q$-Bernoulli numbers $\beta
_{n}(q)$\ and $F_{q}(x,t)$ of $\beta _{n}(x:q)$\ ( $n\geq 0$ ),
respectively, are given as follows : 
\begin{equation}
F_{q}(t)=\frac{q-1}{\log q}\exp (\frac{t}{1-q})-t\sum_{n=0}^{\infty
}q^{n}e^{[n]t}=\sum_{n=0}^{\infty }\frac{\beta _{n}(q)t^{n}}{n!},
\label{EQ-2.1}
\end{equation}%
By using (\ref{q-2}) and (\ref{EQ-2.1}), we have

\begin{equation*}
\lim_{q\rightarrow 1}\beta _{n}(q)=B_{n},\text{ and }\lim_{q\rightarrow
1}F_{q}(t)=F(t).
\end{equation*}%
The generating function $F_{q}(x,t)$ of the $q$-Bernoulli polynomials $\beta
_{n}(x:q)$ ( $n\geq 0$ ) is defined analogously as follows: 
\begin{equation}
F_{q}(x,t)=\frac{q-1}{\log q}\exp (\frac{t}{1-q})-t\sum_{n=0}^{\infty
}q^{n+x}e^{[n+x]t}=\sum_{n=0}^{\infty }\frac{\beta _{n}(x:q)t^{n}}{n!}.
\label{EQ-2.2}
\end{equation}%
By using (\ref{eq-1}) and (\ref{EQ-2.2}), we have

\begin{equation*}
\lim_{q\rightarrow 1}\beta _{n}(x:q)=B_{n}(x),\text{ and }\lim_{q\rightarrow
1}F_{q}(x,t)=F(x,t).
\end{equation*}

The series on the right-hand side of (\ref{EQ-2.1}) and (\ref{EQ-2.2}) are
uniformly convergent in the wider sense. Consequently, we shall explicitly
determine the $q$-Bernoulli numbers as follows: 
\begin{equation*}
\beta _{0}(q)=\frac{q-1}{\log q},\text{ \ }q(q\beta (q)+1)^{n}-\beta
_{n}(q)=\left\{ 
\begin{array}{c}
1,\text{ if }n=1 \\ 
0,\text{ if }n>1,%
\end{array}%
\right.
\end{equation*}%
with the usual convention about replacing $\beta ^{n}$ by $\beta _{n}$.

Let $\chi $ be a Dirichlet character of conductor $f\in \mathbb{Z}^{+}$. The
generating function of generalized $q$-Bernoulli numbers attached to $\chi $
is given as follows ( for detail see \cite{Kim4}, \cite{Kim5}, \cite{Kim7}, %
\cite{Kim9}, \cite{Kim15}, \cite{Kim-simsek-sri}, \cite{Kim-simsek-sri}): 
\begin{eqnarray}
F_{q,\chi }(t) &=&-t\sum_{a=1}^{f}\chi (a)\sum_{n=0}^{\infty
}q^{fn+a}e^{[fn+a]t}  \label{EQ-2.3} \\
&=&-t\sum_{n=0}^{\infty }\chi (n)q^{n}e^{[n]t}  \notag \\
&=&\sum_{n=0}^{\infty }\beta _{n,\chi }(q)\frac{t^{n}}{n!}.  \notag
\end{eqnarray}%
where the coefficients, $\beta _{n,\chi }(q)$ ( $n\geq 0$ ) are called
generalized $q$-Bernoulli numbers with a Dirichlet character. We note from
the definitions in (\ref{eq-3}) and (\ref{EQ-2.3}) that

\begin{equation*}
\lim_{q\rightarrow 1}\beta _{n,\chi }(q)=B_{n,\chi },
\end{equation*}%
and 
\begin{equation*}
\lim_{q\rightarrow 1}F_{q,\chi }(t)=F_{\chi }(t)=\sum_{a=1}^{f}\chi (a)\frac{%
te^{at}}{e^{tf}-1}=\sum_{n=0}^{\infty }B_{n,\chi }\frac{t^{n}}{n!}.
\end{equation*}

Generating function of generalized $q$-Bernoulli polynomials is associated
with a Dirichlet character as follows( \cite{Kim4}, \cite{Kim5}, \cite{Kim7}%
, \cite{Kim9}, \cite{Kim15}, \cite{Kim-simsek-sri}): 
\begin{eqnarray}
F_{q,\chi }(x,t) &=&q^{x}te^{-[x]t}\sum_{n=0}^{\infty }\chi
(n)q^{n}e^{[n]q^{x}t}  \label{EQ-2.6} \\
&=&-t\sum_{n=0}^{\infty }\chi (n)q^{n+x}e^{[n+x]t}  \notag \\
&=&\sum_{n=0}^{\infty }\beta _{n,\chi }(x:q)\frac{t^{n}}{n!}.  \notag
\end{eqnarray}%
By using (\ref{EQ-2.2}), (\ref{EQ-2.3}) and (\ref{EQ-2.6}), we obtain 
\begin{eqnarray*}
F_{q,\chi }(x,t) &=&-t\sum_{n=0}^{\infty }\chi (n)q^{n+x}e^{[n+x]t} \\
&=&-e^{[x]t}q^{x}t\sum_{a=1}^{f}\chi (a)\sum_{n=0}^{\infty
}q^{fn+a}e^{[fn+a]tq^{x}}.
\end{eqnarray*}%
Let $\chi $ be a Dirichlet character of conductor $f\in \mathbb{Z}^{+}$.
Then 
\begin{equation*}
\beta _{n,\chi }(x:q)=\frac{1}{[f]^{1-n}}\sum_{a=0}^{f-1}\chi (a)\beta _{n}(%
\frac{a+x}{f}:q^{f})
\end{equation*}%
(For detail see \cite{Kim5}, \cite{Kim7}, \cite{Kim8}, \cite{Kim9}, \cite%
{Kim13}, \cite{Kim14}, \cite{Kim15}, \cite{Kim16}, \cite{Kim18}, \cite%
{Kim-simsek-sri}).

In \cite{Kim-simsek-sri}, T. Kim, Y. Simsek and H. M. Srivastava gave new
generating functions which produce new definitions of the Barnes-type
Changhee $q$-Bernoulli polynomials and the generalized Barnes-type Changhee $%
q$-Bernoulli numbers with attached to Dirichlet character. These generating
functions are very important in case of multiple zeta function. Therefore,
by using these generating functions, they proved relation between the
Barnes-type Changhee $q$-zeta function and the Barnes-type Changhee $q$%
-Bernoulli numbers.

Let $w,w_{1},w_{2},\ldots ,w_{r}$ be complex numbers such that $w_{i}\neq 0$
for $i=1,2,\ldots ,r$. \ Kim (\cite{Kim4}, \cite{Kim7}, \cite{Kim13}) and
Kim, Simsek and Srivastava \cite{Kim-simsek-sri} defined\ the\ Barnes-type
of Changhee $q$-Bernoulli polynomials of $w$ with parameters $w_{1}$ as
follows: 
\begin{eqnarray}
F_{q}(w,t &\mid &w_{1})=\frac{q-1}{\log q}e^{\frac{t}{1-q}%
}-w_{1}t\sum_{n=0}^{\infty }q^{w_{1}n+w}e^{[w_{1}n+w]t}  \label{EQ-3.1} \\
&=&\sum_{n=0}^{\infty }\frac{\beta _{n}(w:q\mid w_{1})t^{n}}{n!}\text{ \ (}%
\mid t\mid <2\pi \text{ ),}  \notag
\end{eqnarray}%
where the coefficients, $\beta _{n}(w:q\mid w_{1})$ ( $n\geq 0$ ) are called
Barnes-type of Changhee $q$-Bernoulli polynomials in $w$ with parameters $%
w_{1}$.

We note that

\begin{equation*}
\lim_{q\rightarrow 1}\beta _{n}(w:q\mid w_{1})=w_{1}^{n}\beta _{n}(w),\text{
and }\lim_{q\rightarrow 1}F_{q}(w,t\mid w_{1})=\frac{w_{1}t}{e^{w_{1}t}-1}%
e^{wt},
\end{equation*}%
where $\beta _{n}(w)$ are the ordinary Barnes Bernoulli polynomials.

By using (\ref{EQ-3.1}), we easily obtain\cite{Kim4}, \cite{Kim7}, \cite%
{Kim13}, \cite{Kim-simsek-sri} 
\begin{equation*}
\beta _{n}(w:q\mid w_{1})=\frac{1}{(1-q)^{n}}\sum_{l=0}^{n}\left( 
\begin{array}{c}
n \\ 
l%
\end{array}%
\right) q^{lw}(-1)^{l}\frac{lw_{1}}{[lw_{1}]}.
\end{equation*}%
If $w=0$ in the above, then%
\begin{equation*}
\beta _{n}(0:q\mid w_{1})=\beta _{n}(q\mid w_{1}),
\end{equation*}%
where $\beta _{n}(q\mid w_{1})$ are called\ the Barnes-type Changhee $q$%
-Bernoulli numbers with parameter $w_{1}$.

Let $\chi $ be the Dirichlet character with conductor $f$. Then the
generalized Barnes-type Changhee $q$-Bernoulli numbers with attached to $%
\chi $ are defined as follows (\cite{Kim4}, \cite{Kim7}, \cite{Kim13}, \cite%
{Kim-simsek-sri}): 
\begin{eqnarray*}
F_{q,\chi }(t &\mid &w_{1})=-w_{1}t\sum_{n=1}^{\infty }\chi
(n)q^{w_{1}n}e^{[w_{1}n]t} \\
&=&\sum_{n=0}^{\infty }\frac{\beta _{n,\chi }(q\mid w_{1})t^{n}}{n!}\text{ \
(}\mid t\mid <2\pi \text{ ).}
\end{eqnarray*}

\begin{theorem}
(\cite{Kim-simsek-sri}) Let $\chi $ be a Dirichlet character of conductor $%
f\in \mathbb{Z}^{+}$. Then 
\begin{equation}
\beta _{n,\chi }(x:q\mid w_{1})=\frac{1}{[f]^{1-n}}\sum_{a=0}^{f-1}\chi
(a)\beta _{n}(\frac{x+aw_{1}}{f}:q^{f}\mid w_{1}).  \label{bb2}
\end{equation}
\end{theorem}

The Barnes-type Changhee $q$-zeta functions are defined as follows:

\begin{definition}
(\cite{Kim-simsek-sri}) For $s\in \mathbb{C}$, we have 
\begin{equation}
\zeta _{q}(s,w\mid w_{1})=-\frac{(1-q)^{s}}{s-1}\frac{1}{\log q}%
+w_{1}\sum_{n=0}^{\infty }\frac{q^{w_{1}n+w}}{[w_{1}n+w]^{s}}.
\label{EQ-3.5}
\end{equation}
\end{definition}

\begin{theorem}
(\cite{Kim-simsek-sri}) If $n\in \mathbb{Z}^{+}$, then 
\begin{equation}
\zeta _{q}(1-n,w\mid w_{1})=-\frac{\beta _{n}(w:q\mid w_{1})}{n}.
\label{bb1}
\end{equation}
\end{theorem}

\begin{remark}
$\zeta _{q}(s,w\mid w_{1})$ is analytic continuation\ in $\mathbb{C}$\ with
only simple pole at $s=1$.
\end{remark}

\section{The two-variable Dirichlet-type Changhee $q$-$L$-function}

Let (\cite{Kim4}, \cite{Kim5}, \cite{Kim7}, \cite{Kim9}, \cite{Kim15}, \cite%
{Kim-simsek-sri})%
\begin{eqnarray}
F_{q,\chi }(t,x &\mid &w_{1})=-tw_{1}e^{[x]t}\sum_{n=1}^{\infty }\chi
(n)q^{w_{1}n}e^{[nw_{1}]q^{nw_{1}}t}  \notag \\
&=&-tw_{1}\sum_{n=1}^{\infty }\chi (n)q^{nw_{1}}e^{[nw_{1}+x]t}  \label{b} \\
&=&\sum_{n=1}^{\infty }\frac{\beta _{n}(x:q\mid w_{1})}{n!}t^{n},  \notag
\end{eqnarray}%
where we use $[x+a]=[x]+q^{x}[a]$ in (\ref{b}). We consider the following
contour integral: 
\begin{eqnarray}
\frac{\Gamma (1-s)e^{-\pi is}}{2\pi i}\oint_{C}t^{s-2}F_{\chi ,q}(-t,x &\mid
&w_{1})dt=\frac{1}{\Gamma (s)}\int_{0}^{\infty }t^{s-2}F_{\chi ,q}(-t,x\mid
w_{1})dt  \notag \\
&=&\frac{w_{1}}{\Gamma (s)}\sum_{n=1}^{\infty }\chi
(n)q^{w_{1}n}\int_{0}^{\infty }t^{s-1}e^{-[x+w_{1}n]t}dt  \notag \\
&=&w_{1}\sum_{n=0}^{\infty }\frac{\chi (n)q^{w_{1}n}}{[x+w_{1}n]^{s}},
\label{EQ-4.1}
\end{eqnarray}%
where $C$ denotes a positively oriented (counter-clockwise) circle of radius 
$R$, centered at the origin. The function 
\begin{equation*}
Y(t,x)=\frac{1}{\Gamma (s)}t^{s-2}F_{\chi ,q}(-t,x\mid w_{1})
\end{equation*}%
has pole $t=0$ inside the contour $C.$ Therefore, if we want to integrate $%
Y(t,x)$ function, then we have modify the contour by indentation at this
point. We take indentation as identical small semicircle, which has radius $%
r $, leaving $t=0$.

Thus we arrive at the two-variable \ Dirichlet-type Changhee $q$-$L$%
-function, which is given in Definition 1:

Let $\chi $ be a Dirichlet character of conductor $f\in \mathbb{Z}^{+}$, 
\begin{equation}
L_{q}(s,x\mid \chi ;w_{1})=w_{1}\sum_{n=0}^{\infty }\frac{\chi (n)q^{w_{1}n}%
}{[x+w_{1}n]^{s}}.  \label{EQ-4.2}
\end{equation}%
If we take $q\rightarrow 1$ and $w=1$ in (\ref{EQ-4.2}), then Definition 1
reduces to Definition 3, that is 
\begin{equation*}
\lim_{q\rightarrow 1}L_{q}(s,x\mid \chi ;1)=\sum_{n=0}^{\infty }\frac{\chi
(n)}{(x+n)^{s}}.
\end{equation*}

The Dirichlet-type Changhee $q$-$L$-function and the Hurwitz-type Changhee $q
$-zeta function are closely related, too. We give proof of this relation
below.

\begin{proof}[Proof of Theorem 1]
By setting $n=a+kf$, where ( $k=0,1,2,...,\infty $ ; $a=1,2,...,f$ ) in (\ref%
{EQ-4.2}), we have 
\begin{eqnarray*}
L_{q}(s,x &\mid &\chi ;w_{1})=w_{1}\sum_{a=1}^{f}\chi (a)\sum_{k=0}^{\infty }%
\frac{q^{(aw_{1}+kfw_{1})}}{[x+aw_{1}+kfw_{1}]^{s}} \\
&=&w_{1}\sum_{a=1}^{f}\chi (a)\sum_{k=0}^{\infty }\frac{q^{f(\frac{aw_{1}}{f}%
+kw_{1})}}{[f]^{s}[\frac{x+aw_{1}}{f}+kw_{1}:q^{f}]^{s}} \\
&=&[f]^{-s}\sum_{a=1}^{f}\chi (a)\left\{ -\frac{(1-q^{f})^{s}}{s-1}\frac{1}{%
\log q^{f}}+w_{1}\sum_{n=0}^{\infty }\frac{q^{f(\frac{aw_{1}}{f}+kw_{1})}}{[%
\frac{x+aw_{1}}{f}+kw_{1}:q^{f}]^{s}}\right\} .
\end{eqnarray*}%
By using (\ref{EQ-3.5}) in the above equation, we easily arrive at the
desired result.
\end{proof}

\begin{proof}[First proof of Theorem 2]
If we take $s\rightarrow 1-n$ in Theorem 1, where $n$ is a positive integer,
then we have 
\begin{equation*}
L_{q}(1-n,x\mid \chi ;w_{1})=[f]^{n-1}\sum_{a=1}^{f}\chi (a)\zeta
_{q^{f}}(1-n,\frac{x+aw_{1}}{f}\mid w_{1}).
\end{equation*}%
By using (\ref{bb1}) in the above, we get%
\begin{equation*}
L_{q}(1-n,x\mid \chi ;w_{1})=-[f]^{n-1}\sum_{a=1}^{f}\chi (a)\frac{\beta
_{n}(\frac{x+aw_{1}}{f}:q\mid w_{1})}{n}.
\end{equation*}%
By using (\ref{bb2}) in the above, we easily arrive at the following Theorem.
\end{proof}

Now we give proof of Theorem 2 as follows:

\begin{proof}[Second proof of Theorem 2]
Proof of this Theorem similar to that of Theorem 8 in \cite{Kim-simsek-sri}.
Let%
\begin{equation}
B(s,x)=\int_{C}z^{s-2}F_{\chi ,q}(-z,x\mid w_{1})dz,  \label{b1}
\end{equation}%
where $C$ is Hankel's contour along the cut joining the points $z=0$ and $%
z=\infty $ on the real axis, which starts from the point at $\infty $,
encircles the origin ( $z=0$ ) once in the positive (counter-clockwise)
direction, and returns to the point at $\infty $ ( see for details, \cite{E.
T. Wittaker and G. N. Watson} p. 245). Here, as usual, we interpret $z^{s}$
to mean $\exp (s\log z)$, where we assume $\log $ to be defined by $\log t$
on the top part of the real axis and by $\log t+2\pi i$ on the bottom part
of the real axis. We thus find from the definition (\ref{b1}) that 
\begin{equation*}
B(s,x)=(e^{2\pi is}-1)\int_{\varepsilon }^{\infty }t^{s-2}F_{\chi
,q}(-t,x\mid w_{1})dt
\end{equation*}%
\begin{equation*}
+\int_{C_{\varepsilon }}z^{s-2}F_{\chi ,q}(-z,x\mid w_{1})dz,
\end{equation*}%
where $C_{\varepsilon }$ denotes a circle of radius $\varepsilon >0$ (and
centred at the origin), which is described in the positive
(counter-clockwise) direction. Assume first that $\func{Re}(s)>1$. Then 
\begin{equation*}
\int_{C_{\varepsilon }}\rightarrow 0\text{ as }\varepsilon \rightarrow 0,
\end{equation*}%
so we have 
\begin{equation*}
B(s,x)=(e^{2\pi is}-1)\int_{0}^{\infty }t^{s-2}F_{\chi ,q}(-t,x\mid w_{1})dt,
\end{equation*}%
by using (\ref{b}) and (\ref{EQ-4.2}) in the above equation, after some
elementary calculations, we obtain%
\begin{equation*}
B(s,x)=(e^{2\pi is}-1)\Gamma (s)L_{q}(s,x\mid \chi ;w_{1}).
\end{equation*}%
Consequently,%
\begin{equation}
L_{q}(s,x\mid \chi ;w_{1})=\frac{B(s,x)}{(e^{2\pi is}-1)\Gamma (s)},
\label{b2}
\end{equation}%
which, by analytic continuation, holds true for all $s\neq 1$. This
evidently provides us with an analytic continuation of $L_{q}(s,x\mid \chi
;w_{1})$.

Let $s\rightarrow 1-n$ in (\ref{b2}), where $n$ is a positive integer. Since%
\begin{equation*}
e^{2\pi is}=e^{2\pi i(1-n)}=1\text{ \ ( }n\in \mathbb{Z}^{+}\text{ ),}
\end{equation*}%
we have%
\begin{eqnarray}
\lim_{s\rightarrow 1-n}\left\{ (e^{2\pi is}-1)\Gamma (s)\right\}
&=&\lim_{s\rightarrow 1-n}\left\{ \frac{(e^{2\pi is}-1)}{\sin (\pi s)}\frac{%
\pi }{\Gamma (1-s)}\right\}  \notag \\
&=&\frac{2\pi i(-1)^{n-1}}{(n-1)!}\text{ \ ( }n\in \mathbb{Z}^{+}\text{ )}
\label{b3}
\end{eqnarray}%
by means of the familiar reflection formula for $\Gamma (s)$. Furthermore,
since the integrand in (\ref{b1}) has simple pole order $n+1$ at $z=0$,
where also find from the definition (\ref{b1}) with $s=1-n$ that 
\begin{eqnarray}
B(1-n,x) &=&\int_{C}z^{-n-1}F_{\chi ,q}(-z,x\mid w_{1})dz  \notag \\
&=&2\pi i\text{Res}_{z=0}\left\{ \text{ }z^{-n-1}F_{\chi ,q}(-z,x\mid
w_{1})\right\}  \notag \\
&=&(2\pi i)\frac{(-1)^{n}}{n!}\beta _{n}(x:q\mid w_{1}),  \label{b4}
\end{eqnarray}%
where we have made of the power-series representation in (\ref{b}). Thus by
Cauchy Residue Theorem, we easily \ arrive at the desired result upon
suitably combining (\ref{b3}) and (\ref{b4}) with (\ref{b2}).
\end{proof}

Now, we define generalized multiple Changhee $q$-Bernoulli numbers attached
to the Drichlet character $\chi $. We also construct the two-variable
Dirichlet-type multiple Changhee $q$-$L$-functions. We then give relation
between the two-variable Dirichlet-type multiple Changhee $q$-$L$-functions
and the generalized multiple Changhee $q$-Bernoulli numbers as well.

The generalized multiple Changhee $q$-Bernoulli numbers attached to the
Drichlet character $\chi $ are defined by means of the following generating
function: 
\begin{eqnarray}
F_{q,\chi }^{(r)}(t,x &\mid &w_{1},...,w_{r})=(-t)^{r}\left(
\prod_{j=1}^{r}w_{j}\right)  \notag \\
&&\times \sum_{n_{1},...,n_{r}=1}^{\infty }\left( \prod_{k=1}^{r}\chi
(n_{k})\right) q^{\left( \sum_{m=1}^{r}w_{m}n_{m}\right)
}e^{[x+\sum_{m=1}^{r}w_{m}n_{m}]t}  \notag \\
&=&\sum_{n=0}^{\infty }B_{n,\chi }^{(r)}(x,q\mid w_{1},...,w_{r})\frac{t^{n}%
}{n!}\text{ \ (}\mid t\mid <2\pi \text{ ),}  \label{EQ-8.1}
\end{eqnarray}%
where $w_{1},...,w_{r}\in \mathbb{R}^{+}\mathbb{,}$ $r\in \mathbb{Z}^{+}$.

Here, we can now construct the two-variable Dirichlet-type multiple Changhee 
$q$-$L$-function. By using the Mellin transformation and Residue Theorem in (%
\ref{EQ-8.1}), we obtain 
\begin{eqnarray}
\frac{1}{\Gamma (s)}\int_{0}^{\infty }t^{s-1-r}F_{q,\chi }^{(r)}(-t,x &\mid
&w_{1},...,w_{r})dt  \notag \\
&=&\left( \prod_{j=1}^{r}w_{j}\right) \sum_{n_{1},...,n_{r}=1}^{\infty
}\left( \prod_{k=1}^{r}\chi (n_{k})\right) q^{\left(
\sum_{m=1}^{r}w_{m}n_{m}\right) }  \label{EQ-8.5} \\
&&\times \frac{1}{\Gamma (s)}\int_{0}^{\infty
}t^{s-1}e^{-[x+\sum_{m=1}^{r}w_{m}n_{m}]t}dt.  \notag
\end{eqnarray}%
By using (\ref{EQ-8.5}), we can arrive at the definition of the two-variable
Dirichlet-type multiple Changhee $q$-$L$-functions, which is given in
Definition 2, we also give this relation as follows:

For a Dirichlet character $\chi $ with conductor $f\in \mathbb{Z}^{+}$, 
\begin{equation}
L_{q,r}(s,\chi \mid w_{1},...,w_{r})=\left( \prod_{j=1}^{r}w_{j}\right)
\sum_{n_{1},n_{2},...,n_{r}=1}^{\infty }\frac{\left( \prod_{k=1}^{r}\chi
(n_{k})\right) q^{\left( \sum_{m=1}^{r}w_{m}n_{m}\right) }}{%
^{[x+\sum_{m=1}^{r}w_{m}n_{m}]^{s}}}.  \label{EQ-8.6}
\end{equation}

By (\ref{EQ-8.5}) and (\ref{EQ-8.6}), the Theorem 3 provides a relationship
between $L_{q,r}(s,\chi )$ and the Barnes-type multiple Changhee $q$-zeta
functions, is defined by 
\begin{equation}
\zeta _{q,r}(s,w\mid
w_{1},w_{2},...,w_{r})=\sum_{n_{1},n_{2},...,n_{r}=0}^{\infty }\frac{%
q^{w+n_{1}w_{1}+n_{2}w_{2}+...+n_{r}w_{r}}}{%
^{[w+n_{1}w_{1}+n_{2}w_{2}+...+n_{r}w_{r}]^{s}}},\text{ }  \label{EQ-7.7}
\end{equation}%
where $\func{Re}(w)>0,$ $q\in C$ \ with $\mid q\mid <1.$

(For the Barnes-type Changhee multiple $q$-zeta functions see\ \cite{Kim7}, %
\cite{Kim15}, \cite{Kim-simsek-sri}).

\begin{proof}[Proof of Theorem 3]
Proof of Theorem 3 runs parallel to that of Theorem 1 above, so we choose to
omit the details involved. By setting $n_{j}=a_{j}+n_{j}f$, ( $j\in \left\{
1,2,...,r\right\} $, $n_{j}=0,1,...,\infty $, and $a_{j}=1,2,...,f$ \ ) in (%
\ref{EQ-7.7}), we easily arrive at the following Theorem.
\end{proof}

By using (\ref{EQ-8.1}) to Theorem 3, the numbers $L_{q,r}(-n,\chi \mid
w_{1},...,w_{r})$, ( $n>0$ ) are given explicitly by Theorem 4 below.

\begin{proof}[Proof of Theorem 4]
Proof of Theorem 4 runs parallel to that of Theorem 2 above, so we choose to
omit the details involved. If we take $s\rightarrow 1-n$ in (\ref{EQ-8.5}),
where $n$ is a positive integer, then by using the Mellin transformation and
Residue Theorem in (\ref{EQ-8.1}), we easily arrive at the desired result.
\end{proof}

\begin{acknowledgement}
The first author wishes to thank Professor Taekyun Kim and Professor
Seog-Hoon Rim for their hospitality, valuable cowork seminars and financial
support when he was staying at the Kongju National University in 01-07
-2004, 14-08-2004, S. Korea and Kyungpook National University S. Korea in
23-01-2005, 07-02-2005, respectively. He also wishes to thank Mathematicians
and professional colleagues of the Kongju National University and the
Kyungpook National University for their hospitality and friendship.
\end{acknowledgement}


\begin{thebibliography}{99}
\bibitem{Apostol} T. M. Apostol, \textit{Introduction to analytic number
theory}, Springer-Verlag, New York, 1976.

\bibitem{Barnes} E. W. Barnes, On theory of the multiple gamma functions, 
\textit{Trans. Camb. Philos. Soc.}, \textbf{19} (1904), 374-425.

\bibitem{Carlitz} L. Carlitz, $q$-Bernoulli numbers and polynomials, \textit{%
Duke Math.}, \textbf{15} (1948), 987-1000.

\bibitem{Cherednik} I. Cherednik, On $q$-analogues of the Riemann's zeta
function, \textit{Selecta Math.}, \textbf{7} (2001), 447-491.

\bibitem{E. Friedman and S. Ruijsenaars} E. Friedman and S. Ruijsenaars,
Shintani-Barnes zeta and gamma functions, \textit{Adv. in Math.}, \textbf{187%
} (2004), 362-395.

\bibitem{Kim4} T. Kim, An invariant $p$-adic integral associated with Daehee
Numbers, \textit{Integral Transform. Spec. Funct.}, \textbf{13} (2002),
65-69.

\bibitem{Kim5} T. Kim, $q$-Volkenborn integration, \textit{Russ. J. Math
Phys.}, \textbf{19} (2002), 288-299.

\bibitem{Kim7} T. Kim, Non-archimedean $q$-integrals associated with
multiple Changhee $q$-Bernoulli Polynomials, \textit{Russ. J. Math Phys.}, 
\textbf{10} (2003), 91-98.

\bibitem{Kim8} T. Kim, $q$-Riemann zeta function, \textit{Internat. J. Math.
Sci.}, \textbf{2004} (2003), 185--192.

\bibitem{Kim9} T. Kim, On Euller-Barnes multiple zeta functions, \textit{%
Russ. J. Math Phys}., \textbf{10} (2003), 261-267.

\bibitem{Kim13} T. Kim, A note on the $q$-multiple zeta function, \textit{%
Adv. Stud. Contep. Math.}, \textbf{8} (2004), 111-113.

\bibitem{Kim14} T. Kim, $p$-adic $q$-integrals associated with the
Changhee-Barnes' $q$-Bernoulli Polynomials, \textit{Integral Transform.
Spec. Funct.}, \textbf{15} (2004), 415-420.

\bibitem{Kim15} T. Kim, Analytic continuation of multiple $q$-zeta functions
and their values at negative integers, \textit{Russ. J. Math Phys.}, \textbf{%
11} (2004), 71-76.

\bibitem{Kim16} T. Kim, A new approach to $q$-zeta function, \textit{%
arXiv:math.NT/\textbf{0502005} v1 1 Feb.} (2005).

\bibitem{Kim18} T. Kim, A note on multiple zeta function, JP. J. Algebra
Number Theory and Appl. 3 (2003), 471-476.

\bibitem{Kim-Rim1} T. Kim and S. -H. Rim, A note on two variable Dirichlet $%
L $-function, \textit{Adv. Stud. Contep. Math.}, \textbf{10} (2005), 1-7.

\bibitem{Kim -simsek } T. Kim and Y. Simsek, Analytic continuation of the
multiple Daehee $q$-$l$-functions associated with Daehee numbers, \textit{to
appear Russian J. Math}. (2005).

\bibitem{Kim-simsek2} T. Kim Power series and asymptotic series associated
with the $q$-analogue \ of two-variable $p$-adic $L$-function, \textit{Russ.
J. Math Phys.}, \textbf{12 (2) }(2005).

\bibitem{Kim-simsek-sri} T. Kim, Y. Simsek and H. M. Srivastava, $q$%
-Bernoulli numbers and polynomials associated with multiple $q$-zeta
functions and basic $L$-series, \textit{arXiv:math.NT/\textbf{0502019} v%
\textbf{1} 10 Feb.} (2005).

\bibitem{K. Matsumoto} K. Matsumoto, The analytic continuation and the
asyptotic behaviour of certain multiple zeta-function I, \textit{J. Number
Theory,} \textbf{101} (2003), 223-243.

\bibitem{C. A. Nelson and M. G. Gartley} C. A. Nelson and M. G. Gartley, On
the zeros of the $q$-analogue exponential function, \textit{J. Phys. A: Gen.
Math}., \textbf{27} (1994), 3857-3881.

\bibitem{C. A. Nelson and M. G. Gartley-1} C. A. Nelson and M. G. Gartley,
On the two $q$-analogue logaritmic functions: $\ln _{q}(w),\ln (\ln _{q}(z))$%
, \textit{J. Phys. A: Gen. Math.}, \textbf{24} (1996), 8099-8115.

\bibitem{M. Nishizawa} M. Nishizawa, On a $q$-analogue of the multiple gamma
functions, \textit{Lett. Math. Phys.}, \textbf{37} (1996), 2001-2009.

\bibitem{K. Ota} K. Ota, On Kummer-type congruences for derivatives of
Barnes' multiple Bernoulli Polynomials,\textit{\ J. Number Theory,}\textbf{\
92 }(2002), 1-36.

\bibitem{T. M. Rassia and H. M. Srivastava} T. M. Rassias and H. M.
Srivastava, Some classes of infinite series associated with the Riemann zeta
and polygamma functions and generalized harmonic numbers, \textit{Appl.
Math. Computation}, \textbf{131} (2002), 593-605.

\bibitem{S. N. M. Ruijsenaars} S. N. M. Ruijsenaars, On Barnes' multiple
zeta function and gamma functions, \textit{Adv. in Math.}, \textbf{156}
(2000), 107-132.

\bibitem{K. Shiratani and S. Yamamoto} K. Shiratani and S. Yamamoto, On a $p$%
-adic interpolation function for the Euler numbers and its derivative, 
\textit{Mem. Fac. Kyushu Uni.}, \textbf{39} (1985), 113-125.

\bibitem{Y. Simsek} Y. Simsek, Theorems on twisted $L$-functions and twisted
Bernoulli numbers, to appear \textit{Proc. Jangjeon Math. Soc.}

\bibitem{H. M. Srivastava and P. W. Karlsson} H. M. Srivastava and P. W.
Karlsson, \textit{Multiple Gaussian Hypergeometric Series}, Halsted Press
(Ellis Horwood Limited, Chichester), John Wily and Sons, New York,
Chichester, Brisbane and Toronto, 1985.

\bibitem{P. T. Young} P. T. Young, On the behavior of some two-variable $p$%
-adic $L$-function, \textit{J. Number Theory}, \textbf{98} (2003), 67-86.

\bibitem{Washington} L. C. Washington, \textit{Introduction to Cyclotomic
Fields}, Springer-Verlag and New York, 1997.

\bibitem{E. T. Wittaker and G. N. Watson} E. T. Whittaker and G. N. Watson, 
\textit{A course of modern Analysis}, Cambridge Univ. Press, London and New
York, 1927.
\end{thebibliography}
\end{document}